\documentclass[10pt,a4]{article}
\usepackage{latexsym}
\usepackage{amsfonts}
\usepackage[active]{srcltx} 
\usepackage{color}
\usepackage{hyperref}

\usepackage{amsmath}
\usepackage{amssymb}

\newcommand{\ba}{\begin{array}}\newcommand{\ea}{\end{array}}

\newcommand{\ns}{\rm}

\newcommand{\nse}{\kern-3pt\ns$=$}\newcommand{\qd}{\hfill$\Box$\medbreak}

\newcommand{\ext}{\raise1pt\hbox{$\ts\bigwedge$}}

\newcommand{\ts}{\textstyle}
\newcommand{\rf}[1]{(\ref{#1})}
\newcommand{\chii}{\raise2pt\hbox{$\chi$}}

\newcommand{\Fg}{\mbox{${\cal F}\kern-2pt_g$}}

\newcommand{\Mg}{\mbox{${\cal M}\kern-2pt_g$}}
\newcommand{\Ng}{\mbox{${\cal N}\kern-2pt_g$}}
\newcommand{\V}{V\kern-1pt}

\newcommand{\Gg}{\mbox{${\cal G}\kern-2pt_g$}}
\newcommand{\cir}{\raise1.6pt\hbox{\footnotesize$\circ$}}
\newcommand{\tr}{\mbox{\bf tr}}

\newcommand{\diver}{\mbox{\ns div}}
\newcommand{\grad}{\mbox{\ns grad}}

\newcommand{\Res}[2]{\hbox{\ns Res}\kern-16pt\lower5pt\hbox{\footnotesize$_{#1}$}\kern2pt\left[#2\right]}
\newcommand{\qk}{quaternion-K\"ahler\kern2pt}\renewcommand{\,}{\kern1pt}

\newcommand{\End}{{\rm End}}

\renewcommand{\span}{{\rm span}}

\newcommand{\dirac}{/\kern-6pt\partial}

\newcommand{\lra}{\longrightarrow}
\renewcommand{\ts}{\textstyle}
0
\newtheorem{theo}{Theorem}[section]
\newtheorem{defi}{Definition}[section]\newtheorem{lemma}{Lemma}[section]

\newtheorem{corol}{Corollary}[section]
\newtheorem{prop}{Proposition}[section]\def\frac#1#2{{#1\over#2}}
\def\be#1\ee{\begin{equation}#1\end{equation}}

\begin{document}
\title{Spinorially twisted Spin structures, I: curvature identities and eigenvalue
estimates
}

\author{
Malors Espinosa\footnote{Universidad de Guanajuato, Guanajuato, C.P. 36240 Mexico}
\,\,\, and
Rafael Herrera\footnote{Centro de
Investigaci\'on en Matem\'aticas, A. P. 402,
Guanajuato, Gto., C.P. 36000, M\'exico. E-mail: rherrera@cimat.mx}
\footnote{Partially supported by 
grants of CONACyT and LAISLA (CONACyT-CNRS)}}

\date{\today}

\maketitle

{
\abstract{

We define (higher rank) spinorially twisted spin structures 
and deduce various curvature identities as well as estimates for the
eigenvalues 
of the corresponding twisted Dirac operators. 

}
}

\allowdisplaybreaks

\section{Introduction}

The aim of this note is to introduce 
(higher rank) spinorially twisted Spin structures and prove various curvature formulas and
eigenvalue estimates for the corresponding twisted Dirac operators. Such formulas and estimates are the higher
rank analogues of those proved by Hitchin \cite{Hitchin}, Friedrich \cite{friedrich}, Hijazi \cite{Hijazi}
and many others.
We begin by noticing that a Spin$^c$ structure on a Riemannian $n$-dimensional
manifold $M$ consists of the coupling of a (locally
defined) $Spin(n)$-structure and an auxiliary (locally defined)
$U(1)=Spin(2)$ structure \cite{friedrich} and,
similarly, a Spin$^q$ structure on $M$ consists of the coupling of a
(locally defined) $Spin(n)$-structure and an auxiliary (locally defined)
$Sp(1)=Spin(3)$ structure \cite{Nagase}.
Here, we consider analogous twistings with other $Spin(r)$ groups, $r\geq2$, 
in an attempt to develop spinorial techniques to study the geometry of
manifolds which are neither Spin, nor Spin$^c$, nor Spin$^q$.
Although the idea of ``twisting'' is a classical one, here we will take advantage of the spin
geometry (and the Clifford algebra representation) carried by the spinorial twists.

A twisted Spin group is defined as
\[Spin^r(n)= Spin(n)\times_{\mathbb{Z}_2} Spin(r),\]
where $n,r\in\mathbb{N}$. Each factor has a standard spin representation $\Delta_n$ and $\Delta_r$
respectively, so that we can consider the representation of $Spin^r(n)$
\[\Delta_n\otimes \Delta_r,\]
and more generally
\[\Delta_n\otimes \Delta_r^{\otimes m}\]
for odd $m\in\mathbb{N}$.
These are the representations that we will use to associate twisted spinor vector bundles to 
an $n$-dimensional oriented Riemannian manifold $M$ admitting a $Spin^r(n)$
structure (cf. Definition \ref{defi-twisted-spin-structure}), whose sections will be
called spinor
fields, or simply spinors. By choosing a connection on the induced/auxiliary $SO(r)$ principal
bundle on $M$, together with the Levi-Civita connection on $TM$, one can construct a connection on
the twisted spinor bundles and talk about parallel and Killing spinors, as well as twisted
Dirac operators and connection Laplacians.

As in the Spin, Spin$^c$ and Spin$^q$ cases, one can prove (spinorial) curvature identities
(cf. \rf{curvature-twisted-spin}, \rf{Ricci-curvature-identity},
\rf{scalar-curvature-identity}). Such identities render new formulas for the Ricci and the
scalar curvatures of $M$ in the presence of either a parallel spinor (cf. Theorem
\ref{Ricci-parallel-spinor}) or a
Killing spinor (cf. Theorem \ref{Ricci-Killing-spinor}), thanks to the introduction of 
local 2-forms associated to any spinor (cf. Definition \ref{defi-2-forms}).
These formulas involve not only the curvature of the
connection on the auxiliary bundle, but also the aforementioned 2-forms, whose appearance is made
possible by the 
introduction of the spinorial twist. Clearly, the existence of parallel or Killing
spinors on a manifold will impose restrictions on its geometry (holonomy), but this topic will be treated
elsewhere
\cite{Herrera-Santana}.

We also present a twisted version of the Schr\"odinger-Lichnerowicz formula involving the twisted
Dirac operator and the connection Laplacian, and use the Bochner technique to derive
various eigenvalue estimates analogous to those contained in \cite{Hitchin,friedrich,Hijazi,Nakad}.
We would like to
point out that our estimates involve the tensor power $m$ of the twisting bundle (see Corollary 
\ref{corol-no-harmonic-spinors}). 
On the other hand, Corollaries \ref{not-massless-Dirac-spinor} and
\ref{Killing-lower-bound} give new spinor specific criteria to check whether a spinor is parallel, or whether
a number can be a Killing constant in
terms of the associated 2-forms of the spinor and the curvature of the auxiliary bundle.

The paper is organized as follows.
In Section \ref{preliminaries} we recall basic material \cite{friedrich} regarding
Clifford algebras, Clifford multiplication, Spin groups, etc., and define the twisted Spin
groups, twisted spin representations and the 2-forms associated to spinors. We also recall some
material about connections on (twisted) spin bundles and define the twisted Dirac operator and
connection Laplacian.
In Section \ref{results}, we prove various spinorial curvature identities analogous to those in the
Spin$^c$ case (Subsection \ref{curvature-calculations}) and apply them in the cases of parallel and
Killing spinors (Subsections \ref{section-parallel-spinors} and \ref{section-killing-spinors}). We develop in
detail a particular example of
parallel twisted spinor on any oriented Riemannian manifold (cf. Proposition \ref{universal-parallel-spinor}),
with the interesting feature
that its associated 2-forms generate a
copy of $\mathfrak{so}(n)$ (cf. Proposition \ref{second-outstanding-feature}). Finally. we prove the
twisted Schr\"odinger-Lichnerowicz formula (Theorem \ref{theo-SL-formula}) and prove several
corollaries (Subsection
\ref{section-twisted-SL-formula}).

{\em Acknowledgements}. The second named author would like to thank the International
Centre for Theoretical Physics (Italy) and the Institut des Hautes \'Etudes Scientifiques (France)
for their hospitality and support.

\section{Preliminaries}\label{preliminaries}

In this section, we recall basic material that can be consulted in \cite{friedrich} and 
define the various twisted objects that will be used throughout.

\subsection{Clifford algebra, spin group and representation}

Let $Cl_n$ denote the $2^n$-dimensional real Clifford algebra generated by the orthonormal vectors
$e_1, e_2, \ldots, e_n\in \mathbb{R}^n$ subject to
the relations
$e_i e_j + e_j e_i = -2\delta_{ij}$,
and
$\mathbb{C}l_n=Cl_n\otimes_{\mathbb{R}}\mathbb{C}$
denote its complexification. Recall that
\[\mathbb{C}l_n\cong \left\{
                     \begin{array}{ll}
                     \End(\mathbb{C}^{2^k}), & \mbox{if $n=2k$,}\\
                     \End(\mathbb{C}^{2^k})\oplus\End(\mathbb{C}^{2^k}), & \mbox{if $n=2k+1$,}
                     \end{array}
\right.
\]
and that the space of (complex) spinors is defined as
\[\Delta_n:=\mathbb{C}^{2^k}=\underbrace{\mathbb{C}^2\otimes \ldots \otimes \mathbb{C}^2}_{k\,\,\,\,times}.\]
The map
\[\kappa:\mathbb{C}l_n \lra \End(\mathbb{C}^{2^k})\]
is defined to be either the aforementioned isomorphism for $n$ even, or the isomorphism followed
by the projection onto the first summand for $n$ odd.
In order to make $\kappa$ explicit, consider the following matrices
\[Id = \left(\begin{array}{ll}
1 & 0\\
0 & 1
      \end{array}\right),\quad
g_1 = \left(\begin{array}{ll}
i & 0\\
0 & -i
      \end{array}\right),\quad
g_2 = \left(\begin{array}{ll}
0 & i\\
i & 0
      \end{array}\right),\quad
T = \left(\begin{array}{ll}
0 & -i\\
i & 0
      \end{array}\right).
\]
In terms of the generators $e_1, \ldots, e_n$ of the Clifford algebra, $\kappa$ can be
described explicitly as follows,
\begin{eqnarray}
e_1&\mapsto& Id\otimes Id\otimes \ldots\otimes Id\otimes Id\otimes g_1,\nonumber\\
e_2&\mapsto& Id\otimes Id\otimes \ldots\otimes Id\otimes Id\otimes g_2,\nonumber\\
e_3&\mapsto& Id\otimes Id\otimes \ldots\otimes Id\otimes g_1\otimes T,\nonumber\\
e_4&\mapsto& Id\otimes Id\otimes \ldots\otimes Id\otimes g_2\otimes T,\nonumber\\
\vdots && \dots\nonumber\\
e_{2k-1}&\mapsto& g_1\otimes T\otimes \ldots\otimes T\otimes T\otimes T,\nonumber\\
e_{2k}&\mapsto& g_2\otimes T\otimes\ldots\otimes T\otimes T\otimes T,\nonumber
\end{eqnarray}
and, if $n=2k+1$, 
\[ e_{2k+1}\mapsto i\,\, T\otimes T\otimes\ldots\otimes T\otimes T\otimes T.\]
The vectors 
\[u_{+1}={1\over \sqrt{2}}(1,-i)\quad\quad\mbox{and}\quad\quad u_{-1}={1\over \sqrt{2}}(1,i),\]
form a unitary basis of $\mathbb{C}^2$ with respect to the standard Hermitian product.
Thus, 
\[\{u_{\varepsilon_1,\ldots,\varepsilon_k}=u_{\varepsilon_1}\otimes\ldots\otimes
u_{\varepsilon_k}\,\,|\,\, \varepsilon_j=\pm 1,
j=1,\ldots,k\},\]
is a unitary basis of $\Delta_n=\mathbb{C}^{2^k}$
with respect to the naturally induced Hermitian product.
We will denote inner and Hermitian products (as well as Riemannian and Hermitian
metrics) by the same symbol $\left<\cdot,\cdot\right>$ trusting that the context will make clear
which product is being used.

By means of $\kappa$ we have Clifford multiplication
\begin{eqnarray*}
\mu_n:\mathbb{R}^n\otimes \Delta_n &\lra&\Delta_n\\ 
x \otimes \phi &\mapsto& \mu_n(x\otimes \phi)=x\cdot\phi :=\kappa(x)(\phi) .
\end{eqnarray*}
A quaternionic structure $\alpha$ on $\mathbb{C}^2$ is given by
\[\alpha\left(\begin{array}{c}
z_1\\
z_2
              \end{array}
\right) = \left(\begin{array}{c}
-\overline{z}_2\\
\overline{z}_1
              \end{array}\right),\]
and a real structure $\beta$ on $\mathbb{C}^2$ is given by
\[\beta\left(\begin{array}{c}
z_1\\
z_2
              \end{array}
\right) = \left(\begin{array}{c}
\overline{z}_1\\
\overline{z}_2
              \end{array}\right).\]
Following \cite[p. 31]{friedrich}, the real and quaternionic structures $\gamma_n$  on
$\Delta_n=(\mathbb{C}^2)^{\otimes
[n/2]}$ are built as follows
\[
\begin{array}{cclll}
 \gamma_n &=& (\alpha\otimes\beta)^{\otimes 2k} &\mbox{if $n=8k,8k+1$}& \mbox{(real),} \\
 \gamma_n &=& \alpha\otimes(\beta\otimes\alpha)^{\otimes 2k} &\mbox{if $n=8k+2,8k+3$}&
\mbox{(quaternionic),} \\
 \gamma_n &=& (\alpha\otimes\beta)^{\otimes 2k+1} &\mbox{if $n=8k+4,8k+5$}&\mbox{(quaternionic),} \\
 \gamma_n &=& \alpha\otimes(\beta\otimes\alpha)^{\otimes 2k+1} &\mbox{if $n=8k+6,8k+7$}&\mbox{(real).}
\end{array}
\]

The Spin group $Spin(n)\subset Cl_n$ is the subset 
\[Spin(n) =\{x_1x_2\cdots x_{2l-1}x_{2l}\,\,|\,\,x_j\in\mathbb{R}^n, \,\,
|x_j|=1,\,\,l\in\mathbb{N}\},\]
endowed with the product of the Clifford algebra.
It is a Lie group and its Lie algebra is
\[\mathfrak{spin}(n)=\mbox{span}\{e_ie_j\,\,|\,\,1\leq i< j \leq n\}.\]
The restriction of $\kappa$ to $Spin(n)$ defines the Lie group representation
\[\kappa_n:=\kappa|_{Spin(n)}:Spin(n)\lra GL(\Delta_n),\]
which is, in fact, special unitary. 
We have the corresponding Lie algebra representation
\[\kappa_{n_*}:\mathfrak{spin}(n)\lra \mathfrak{gl}(\Delta_n).\] 
Both representations can be extended to tensor powers $\Delta_n^{\otimes m}$, $m\in \mathbb{N}$, in
the usual way.

Recall that the Spin group $Spin(n)$ is the universal double cover of $SO(n)$, $n\ge 3$. For $n=2$
we consider $Spin(2)$ to be the connected double cover of $SO(2)$.
The covering map will be denoted by 
\[\lambda_n:Spin(n)\rightarrow SO(n)\subset GL(\mathbb{R}^n).\]
Its differential is given
by $\lambda_{n_*}(e_ie_j) = 2E_{ij}$, where $E_{ij}=e_i^*\otimes e_j - e_j^*\otimes e_i$ is the
standard basis of the skew-symmetric matrices, and $e^*$ denotes the metric dual of the vector $e$.
Furthermore, we will abuse the notation and also denote by $\lambda_n$ the induced representation on
the exterior algebra $\ext^*\mathbb{R}^n$. 

The Clifford multiplication $\mu_n$ 
 is skew-symmetric with respect to the Hermitian product
\begin{equation}
\left<x\cdot\phi_1 , \phi_2\right> =\left<\mu_n(x\otimes \phi_1) , \phi_2\right> 
=-\left<\phi_1 , \mu_n(x\otimes \phi_2)\right>
=-\left<\phi_1 , x\cdot \phi_2\right>, \label{clifford-skew-symmetric}
\end{equation} 
 is  $Spin(n)$-equivariant and can be extended to a $Spin(n)$-equivariant map 
\begin{eqnarray*}
\mu_n:\ext^*(\mathbb{R}^n)\otimes \Delta_n &\lra&\Delta_n\\ 
\omega \otimes \psi &\mapsto& \omega\cdot\psi.
\end{eqnarray*}

\subsection{Spinorially twisted spin groups and representations}

We define the twisted Spin group
$Spin^r(n)$ as follows 
\[
Spin^{r}(n) =  (Spin(n) \times Spin(r))/\{\pm (1,1)\} =
Spin(n) \times_{\mathbb{Z}_2} Spin(r), 
\]
where $r\in\mathbb{N}$ and $r\geq 2$,
which fit into exact sequences
\[1\lra \mathbb{Z}_2\lra Spin^r(n)\xrightarrow{\lambda_n\times\lambda_r} SO(n)\times SO(r)\lra 1,\]
where
\begin{eqnarray*}
 \lambda_n\times\lambda_r: Spin^r(n)&\longrightarrow& SO(n)\times SO(r)\\
\,[g,h]  &\mapsto& (\lambda_n(g),\lambda_r(h)).
\end{eqnarray*}
We will call $r$ the {\em rank} of the twisting.
Note that the groups $Spin^{2}(n)=Spin^c(n)$ and $Spin^{3}(n)=Spin^q(n)$.
The Lie algebra of $Spin^r(n)$  is
\[\mathfrak{spin}^r(n) = \mathfrak{spin}(n) \oplus \mathfrak{spin}(r).\]
Consider the representations
\begin{eqnarray*}
\kappa_n\otimes\kappa_r^m: Spin^r(n)&\longrightarrow& GL(\Delta_n\otimes
\Delta_r^{\otimes m})\\
\,[g,h]  &\mapsto& \kappa_n(g)\otimes\kappa_r^m(h),
\end{eqnarray*}
where $m\in\mathbb{N}$,
which are unitary with respect to the Hermitian metric,
and the map
\begin{eqnarray*}
 \mu_n\otimes\mu_r:\left(\ext^*\mathbb{R}^n\otimes_\mathbb{R} \ext^*\mathbb{R}^r\right)
 \otimes_\mathbb{R} (\Delta_n\otimes \Delta_r) &\longrightarrow& \Delta_n\otimes\Delta_r\\
(w_1 \otimes w_2)\otimes (\psi\otimes \varphi) &\mapsto& 
(w_1\otimes w_2)\cdot (\psi\otimes \varphi)
= (w_1\cdot\psi) \otimes (w_2\cdot \varphi).
\end{eqnarray*}
As in the untwisted case, 
$\mu_n\otimes\mu_r$ is an equivariant homomorphism of $Spin^r(n)$ representations.
Note that we can also take tensor products with more copies of $\Delta_r$ as follows
\begin{eqnarray*}
\mu_r^a:= Id_{\Delta_r}^{\otimes a-1}\otimes \mu_r\otimes
Id_{\Delta_r}^{\otimes m-a}:
\ext^*\mathbb{R}^r
 \otimes_\mathbb{R}  \Delta_r^{m} &\longrightarrow&
\Delta_r^{m}\\
(\beta)\otimes ( \varphi_1\otimes\cdots\otimes
\varphi_{m})
&\mapsto& 
  \varphi_1\otimes\cdots\otimes (\mu_r(\beta\otimes 
\varphi_a))\otimes\cdots\otimes
\varphi_{m},
\end{eqnarray*}
with Clifford multiplication tanking place only in the $a$-th factor. 
We will also write
\[\mu_r^a(\beta\otimes \varphi_1\otimes\cdots\otimes \varphi_{m}) = \mu_r^a(\beta)\cdot
(\varphi_1\otimes\cdots\otimes \varphi_{m}). \]
Notice that
\begin{equation}
\kappa_{r*}^m(f_kf_l) (\varphi_1\otimes\cdots \otimes\varphi_m)=
(\mu_r^1(f_kf_l)\cdot\varphi_1)\otimes\cdots \otimes\varphi_m +\cdots+ \varphi_1\otimes\cdots
\otimes(\mu_r^m(f_kf_l)\cdot\varphi_m).\label{kappa-en-mus} 
\end{equation}

An element $\phi\in\Delta_n\otimes\Delta_r^{\otimes m}$ will be called a {\em
twisted spinor}, or simply a {\em spinor}.

\begin{defi}\label{defi-2-forms}
Let $\phi\in\Delta_n\otimes\Delta_r^{\otimes m}$, $X,Y\in\mathbb{R}^n$ and $(f_1\ldots,f_r)$
an orthonormal frame of $\mathbb{R}^r$.
\begin{itemize}
\item Let
\[\eta_{kl}^{\phi} (X,Y) = {\rm Re}\left< X\wedge Y\cdot \kappa_{r*}^m(f_kf_l)\cdot
\phi,\phi\right>\]
be the real $2$-form associated to the spinor $\phi$ where $1\leq k,l\leq r$.
\item Define the antisymmetric endomorphism
$\hat\eta_{kl}^\phi\in\End^-(\mathbb{R}^n)$ by
\[X\mapsto \hat\eta_{kl}^\phi(X):=(X\lrcorner \,\eta_{kl}^{\phi})^\sharp,\]
where $X\in\mathbb{R}^n$, $1\leq k,l\leq r$, $\lrcorner$ denotes contraction and $^\sharp$ denotes metric
dualization.
\end{itemize}
\end{defi}

In fact, for any $\xi\in\ext^2(\mathbb{R}^n)^*$, we define $\hat\xi\in\End^-(\mathbb{R}^n)$ by
\begin{eqnarray*}
 \hat\xi:\mathbb{R}^n &\lra& \mathbb{R}^n \\
     x &\mapsto& \hat\xi(x) := (x\lrcorner \xi)^\sharp.
\end{eqnarray*}

\begin{lemma}
Let $\phi\in\Delta_n\otimes\Delta_r^{\otimes m}$, $X,Y\in\mathbb{R}^n$, $(f_1\ldots,f_r)$
an orthonormal basis of $\mathbb{R}^r$ and $1\leq k,l\leq r$. Then
\begin{eqnarray}
{\rm Re}\left<  \kappa_{r*}^m(f_kf_l)\cdot
\phi,\phi\right>&=&0,\label{vanishing2}\\
{\rm Re}\left< X\wedge Y\cdot
\phi,\phi\right>&=&0,\label{vanishing3}\\
{\rm Im}\left< X\wedge Y\cdot \kappa_{r*}^m(f_kf_l)\cdot
\phi,\phi\right>&=&0,\label{vanishing4}\\
{\rm Re} \left< X\cdot \phi,Y\cdot\phi \right> 
 &=&   \left<X,Y\right>|\phi|^2, \label{real-part}
\end{eqnarray}
\end{lemma}
{\em Proof}.
By using \rf{clifford-skew-symmetric} repeatedly
\begin{eqnarray*}
 \left<  \mu_{r}^a(f_kf_l)\cdot\phi,\phi\right>
&=& -\overline{\left<   \mu_{r}^a(f_kf_l)\phi,\phi\right>},
\end{eqnarray*}
so that \rf{vanishing2} follows
from \rf{kappa-en-mus}.

For identity \rf{vanishing3}, recall that for $X, Y\in \mathbb{R}^n$ 
\[X\wedge Y =  X\cdot Y + \left<X,Y\right>.\]
Thus
\begin{eqnarray*}
 \left< X\wedge Y\cdot \phi,\phi\right>
 &=& -\overline{\left< X\wedge Y \cdot\phi,\phi\right>}.
\end{eqnarray*}

Identities \rf{vanishing4} and \rf{real-part} follow similarly.
\qd

{\bf Remarks}. 
\begin{itemize}
 \item For $k\not= l$,
\[\eta_{kl}^\phi = - \eta_{lk}^\phi.\]
\item By \rf{vanishing3}, 
\[\eta_{kk}\equiv 0.\]
 \item By \rf{vanishing4}, if $k\not= l$,
\[\eta_{kl}^{\phi} (X,Y) =\left< X\wedge Y\cdot \kappa_{r*}^m(f_kf_l)\cdot
\phi,\phi\right>.\]
\item Note that, depending on the spinor, such 2-forms can actually be identically zero.
\end{itemize}

\begin{lemma}
Any spinor $\phi\in\Delta_n\otimes\Delta_r^{\otimes m}$ defines two maps (extended by linearity)
\begin{eqnarray*}
\ext^2 \mathbb{R}^r&\lra& \ext^2 \mathbb{R}^n\\
f_kf_l &\mapsto& \eta_{kl}^{\phi},
\end{eqnarray*}
and
\begin{eqnarray*}
\ext^2 \mathbb{R}^r&\lra& \End(\mathbb{R}^n)\\
f_kf_l &\mapsto& \hat\eta_{kl}^{\phi}.
\end{eqnarray*}
\end{lemma}
\qd

\subsection{Spinorially twisted spin structures on oriented Riemannian
manifolds}\label{twisted-structures}

\begin{defi}\label{defi-twisted-spin-structure}
Let $M$ be an oriented $n$-dimensional Riemannian manifold, $P_{SO(M)}$ be its principal bundle of
orthonormal frames and $r\in\mathbb{N}$, $r\geq 2$. 
A $Spin^r(n)$ structure on $M$ consists of an auxiliary 
principal $SO(r)$-bundle $P_{SO(r)}$ and
a principal $Spin^r(n)$-bundle $P_{Spin^{r}(n)}$ together with an equivariant $2:1$ covering map
\[\Lambda:P_{Spin^{r}(n)}\lra P_{SO(M)} \tilde{\times} P_{SO(r)},\]
where $\tilde{\times}$ denotes the fibered product,
such that $\Lambda(pg)=\Lambda(p)(\lambda_n\times\lambda_r)(g)$ for all $p\in P_{Spin^{r}(n)}$ and
$g\in Spin^{r}(n)$, where $\lambda_n\times\lambda_r:Spin^r(n)\lra SO(n)\times SO(r)$ denotes the
canonical $2$-fold cover.

A manifold $M$ admitting a $Spin^r(n)$ structure will be called a {\em Spin$^r$ manifold}.
\end{defi}

{\bf Remark}. A Spin$^r$ manifold with trivial $P_{SO(r)}$ auxiliary bundle is a Spin manifold.
Conversely, any Spin manifold admits $Spin^r(n)$ structures with trivial $P_{SO(r)}$ via
the inclusion $Spin(n)\subset Spin^r(n)$ given by the elements $[g,1]$

{\bf Remark}. A Spin$^r$ manifold has the following associated vector bundles:
\begin{eqnarray*}
TM &=& P_{Spin^r(n)}\times_{\lambda_n\times\lambda_r} (\mathbb{R}^n\times\{0\}),\\
F&=&P_{Spin^r(n)}\times_{\lambda_n\times\lambda_r} (\{0\}\times\mathbb{R}^r),\\
S(TM)\otimes S(F)^{\otimes m}&=& P_{Spin^r(n)}\times_{\kappa_n\otimes\kappa_r^m} 
(\Delta_n\otimes\Delta_r^{\otimes m}),
\end{eqnarray*}
where the last bundle is globally defined if
$M$ and $m$ satisfy certain conditions.
Indeed, $S(TM)\otimes S(F)^{\otimes m}$ is defined if one of the following options holds:
\begin{itemize}
 \item $M$ is a non-Spin Spin$^r$ manifold and $m$ is odd. The structure group under
consideration is $Spin^r(n)$.
 \item Both $M$ and $F$ admit Spin structures, and $m\in\mathbb{N}$. The structure group under
consideration is 
$Spin(n)\times Spin(r)$, so that we can consider all representations of the product group.
 \item $M$ is Spin, $F$ is not Spin, and $m$ must be even. In this case, the representation
$\Delta_r^{\otimes m}$ must factor through $SO(r)$ in order to get a globally defined bundle.
Thus, the structure group we need to consider is $Spin(n)\times SO(r)$.
Although this case falls outside the definition of Spin$^r$ structure,
we will still consider it since one can still work with twisted spinors and twisted Dirac operators.
\end{itemize}

{\bf Example}. Let $r = ak + bl$, $a,b\in\mathbb{N}$ and consider the real Grassmannians of oriented subspaces
\[\mathbb{G}r_k(\mathbb{R}^{k+l})={SO(k+l)\over SO(k)\times SO(l)}.\]
There exists a homomorphism
$
SO(k)\times SO(l)\rightarrow Spin^{r}(kl)$ providing a homogeneous
$Spin^{r}(kl)$-structure on the real Grassmannian
$\mathbb{G}r_k(\mathbb{R}^{k+l})$
if 
\begin{eqnarray*}
 a&\equiv& l \pmod 2,\\
b&\equiv& k \pmod 2.
\end{eqnarray*}

\subsection{Covariant derivatives on spinorially twisted Spin bundles}

Let $M$ be a Spin$^r$ $n$-dimensional manifold and $F$ its auxiliary Riemannian vector bundle of
rank $r$.
Assume $F$ is endowed with a covariant derivative $\nabla^F$ (or equivalently, that $P_{SO(F)}$ is
endowed with a connection 1-form $\theta$) and denote by $\nabla$ the
Levi-Civita covariant derivative on $M$. These two derivatives induce the spinor covariant
derivative
\[\nabla^{\theta}: \Gamma(S(TM)\otimes S(F)^{\otimes m})\lra \Gamma(T^*M\otimes
S(TM)\otimes S(F)^{\otimes m})\]
\[\nabla^{\theta} (\psi\otimes \varphi) =d(\psi\otimes\varphi)
+\left[{1\over 2}\sum_{1\leq i<j\leq n}\omega_{ji}\otimes e_ie_j\cdot\psi\right] \otimes\varphi
+\psi\otimes\left[{1\over 2}\sum_{1\leq k<l\leq r}\theta_{kl}\otimes \kappa_{r*}^m(f_kf_l)\cdot
\varphi\right]
,\]
where $\psi\otimes\varphi\in\Gamma(S(TM)\otimes S(F)^{\otimes m})$, 
$(e_1,\ldots,e_n)$ and $(f_1,\ldots,f_r)$ are a local orthonormal frames of $TM$ and $F$
resp., and
$\omega_{ij}$ and $\theta_{kl}$ are the local connection 1-forms for $TM$ (Levi-Civita) and $F$ respectively.

From now on, we will omit the upper and lower bounds on the indices, by declaring $i$ and $j$
to be
the indices for the frame vectors of $M$, and $k$ and $l$ to be the indices for the local frame sections
of $F$. 
Now, for any tangent vectors $X,Y\in T_xM$,
\begin{eqnarray}
R^{\theta}(X,Y)(\psi\otimes\varphi)&=&
\left[{1\over 2} \sum_{i<j} \Omega_{ij}(X,Y) e_ie_j\cdot\psi\right]\otimes\varphi+
\psi\otimes\left[{1\over 2} \sum_{k<l} \Theta_{kl}(X,Y) \kappa_{r*}^m(f_kf_l)\cdot\varphi\right],
\label{curvature-twisted-spin}
\end{eqnarray}
where
\[\Omega_{ij}(X,Y)=\big<R^M(X,Y)(e_i),e_j
\big>\quad\quad\mbox{and}\quad\quad\Theta_{kl}(X,Y)=\big<R^F(X,Y)(f_k),f_l \big>.\]

For $X, Y$ vector fields and $\phi\in\Gamma(S(TM)\otimes S(F)^{\otimes m})$ a spinor field,
we also have the compatibility of the covariant derivative with Clifford multiplication, 
\[
\nabla^{\theta}_X(Y\cdot\phi) = (\nabla_XY)\cdot\phi +
Y\cdot\nabla_X^{\theta}\phi.
\]

\subsection{Twisted differential operators}

In order to simplify notation, let $S=S(TM)\otimes S(F)^{\otimes m}$ and
$\phi\in\Gamma(S)$.

\begin{defi}
The {\em twisted Dirac operator} is the first order differential operator 
$\dirac^\theta=\dirac^{\theta,m}:\Gamma(S)\longrightarrow \Gamma(S)$
defined by
\begin{eqnarray*}
\dirac^\theta(\phi)
&=&\sum_{i=1}^n e_i\cdot\nabla_{e_i}^{\theta}(\phi).
\end{eqnarray*}
\end{defi}
We will generally use the notation $\dirac^\theta$, and 
will use the notation $\dirac^{\theta,m}$ whenever we want to emphasize which tensor power is involved in the
twisted vector bundle being considered. 

{\bf Remark}.
The twisted Dirac operator $\dirac^{\theta}$ is well-defined and formally self-adjoint on compact
manifolds. Moreover, 
if $h\in C^{\infty}(M)$, $\phi\in\Gamma(S)$, we have
\[
\dirac^\theta(h\,\phi) = \grad(h)\cdot \phi +
h\,\dirac^\theta(\phi).
\]
The proofs of these facts are analogous to those for the Spin$^c$ Dirac operator
\cite{friedrich}.

\begin{defi}
The {\em twisted spin connection Laplacian} is the second order differential operator 
$\Delta: \Gamma(S)\rightarrow\Gamma(S)$  defined as
\[
\Delta^{\theta}(\phi) = - \sum_{i = 1}^n
\nabla^{\theta}_{e_i}\nabla^{\theta}_{e_i}(\phi) - \sum_{i = 1}^n
\diver(e_i)\nabla^{\theta}_{e_i}(\phi). 
\]
\end{defi}

\section{Curvature identities, special spinors and twisted Dirac operator's eigenvalue estimates
}\label{results}

Throughout this section, 
let $M$ be a Spin$^r$ $n$-dimensional manifold, $m\in\mathbb{N}$, such that 
$S=S(TM)\otimes S(F)^{\otimes m}$ is globally defined,
$(e_1,...,e_n)$ and $(f_1,...,f_r)$ be local orthonormal frames of $TM$ and $F$ respectively. 

\subsection{Curvature calculations}\label{curvature-calculations}

\begin{prop}
For $X\in\Gamma(TM)$
and $\phi\in\Gamma(S)$, we have
\begin{eqnarray}
\sum_{i = 1}^n e_i\cdot R^{\theta}(X, e_i)(\phi) 
 &=& -\dfrac{1}{2} {\rm Ric}(X)
\cdot \phi + \dfrac{1}{2}\sum_{k<l}(X\lrcorner\Theta_{kl})\cdot
\kappa_{r*}^m(f_kf_l)\cdot\phi.\label{Ricci-curvature-identity}
\end{eqnarray}
\end{prop}
{\em Proof}. 
By \rf{curvature-twisted-spin}, 
if $\phi=\psi\otimes \varphi$,
\[
R^\theta(X,e_\alpha)(\psi\otimes \varphi) = 
\left[{1\over 2} \sum_{i<j} \Omega_{ij}(X,e_\alpha) e_ie_j\cdot\psi\right]\otimes \varphi+
\psi\otimes\left[{1\over 2} \sum_{k<l} \Theta_{kl}(X,e_\alpha) \kappa_{r*}^m(f_kf_l)\cdot
\varphi\right].
\]
Multiply by $e_\alpha$ 
and sum over $\alpha$
\begin{eqnarray*}
\sum_\alpha e_\alpha\cdot R^\theta(X,e_\alpha)(\psi\otimes \varphi) 
&=& 
\left[{1\over 2} \sum_\alpha \sum_{i<j} \Omega_{ij}(X,e_\alpha)  e_\alpha
e_ie_j\cdot\psi\right]\otimes
\varphi+
{1\over 2} \sum_{k<l}\left[\sum_\alpha \Theta_{kl}(X,e_\alpha) e_\alpha \cdot \psi\right]\otimes 
\kappa_{r*}^m(f_kf_l)\cdot
\varphi.
\end{eqnarray*}
The term 
\[{1\over 2} \sum_\alpha \sum_{i<j} \Omega_{ij}(X,e_\alpha)  e_\alpha
e_ie_j  = -{1\over 2} {\rm Ric}(X)
\]
(cf. \cite{friedrich}).
The second term
\begin{eqnarray*}
{1\over 2} \sum_{k<l}\left[\sum_\alpha \Theta_{kl}(X,e_\alpha) e_\alpha \cdot \psi\right]\otimes 
\kappa_{r*}^l(f_kf_l)\cdot\varphi
&=& {1\over 2} \sum_{k<l} (X\lrcorner\Theta_{kl}) \cdot
\kappa_{r*}^m(f_kf_l)\cdot
(\psi\otimes\varphi).
\end{eqnarray*}
\qd

\begin{prop}
Let $\phi\in\Gamma(S)$. Then
\begin{equation}
\sum_{i,j} e_ie_j\cdot R^{\theta}(e_i, e_j)(\phi) 
= {R\over 2}\phi  + \sum_{k<l}\Theta_{kl}\cdot
\kappa_{r*}^m(f_kf_l)\cdot\phi,\label{scalar-curvature-identity}
\end{equation}
where
\[\Theta_{kl}=\sum_{i<j}\Theta_{kl}(e_i,e_j)e_i\wedge e_j.\]
\end{prop}
{\em Proof}.
By \rf{Ricci-curvature-identity},
\[
\sum_{j = 1}^n e_j\cdot R^{\theta}(e_i, e_j)(\phi) = -\dfrac{1}{2} {\rm Ric}(e_i)
\cdot \phi + \dfrac{1}{2}\sum_j\sum_{k<l}\Theta_{kl}(e_i,e_j)e_j\cdot
\kappa_{r*}^m(f_kf_l)\cdot\phi,
\]
By multiplying with $e_i$ and summing over $i$, we get
\begin{eqnarray*}
\sum_{i,j} e_ie_j\cdot R^{\theta}(e_i, e_j)(\phi) 
&=&
 -\dfrac{1}{2} \sum_ie_i\cdot{\rm Ric}(e_i)
\cdot \phi + \dfrac{1}{2}\sum_{k<l}\left[\sum_{i,j}\Theta_{kl}(e_i,e_j)e_ie_j\right]\cdot
\kappa_{r*}^m(f_kf_l)\cdot\phi.
\end{eqnarray*}
Now,
\begin{eqnarray*}
- \sum_i e_i\cdot {\rm Ric}(e_i) 
 &=& {\rm R},
\end{eqnarray*}
where ${\rm R}$ denotes the scalar curvature of $M$, and
for $k$ and $l$ fixed, 
\begin{eqnarray*}
 \sum_{i,j}\Theta_{kl}(e_i,e_j)e_ie_j 
&=& 2\sum_{i<j}\Theta_{kl}(e_i,e_j)e_ie_j
\quad=\quad 2\Theta_{kl}.
\end{eqnarray*}
\qd

Now, let us denote
\begin{eqnarray*}
\Theta&=&\sum_{k<l}\Theta_{kl}\otimes f_kf_l \in \ext^2T^*M\otimes \ext^2 F ,\\ 
\hat{\Theta}&=&\sum_{k<l}\hat{\Theta}_{kl}\otimes f_kf_l \in \End^-(TM)\otimes \ext^2F,
\end{eqnarray*}
and 
\begin{eqnarray*}
\eta^\phi&=&\sum_{k<l}\eta_{kl}^\phi\otimes f_kf_l \in \ext^2T^*M\otimes \ext^2 F ,\\ 
\hat{\eta}^\phi&=&\sum_{k<l}\hat{\eta}_{kl}^\phi\otimes f_kf_l \in \End^-(TM)\otimes \ext^2F.
\end{eqnarray*}
In order to simplify notation, we define
\begin{eqnarray*}
 \left<{\Theta},{\eta}^\phi\right>_0 
 &=&
 \sum_{k<l}\sum_{i<j} \Theta_{kl}(e_i,e_j)\eta_{kl}^\phi(e_i,e_j),\\
 \left<\hat{\Theta},\hat{\eta}^\phi\right>_1 
 &=&
 \sum_{k<l}\tr( \hat{\Theta}_{kl}(\hat{\eta}_{kl}^\phi)^T).
\end{eqnarray*}

\subsection{Parallel spinors}\label{section-parallel-spinors}

\begin{defi}
 A spinor $\phi\in\Gamma(S)$ is said to be parallel if
\[\nabla^\theta_X \phi =0 \]
for all $X\in\Gamma(TM)$.
\end{defi}

\begin{theo}\label{Ricci-parallel-spinor}
 Let $\phi\in\Gamma(S)$ be a non-zero parallel spinor.
Then 
\begin{enumerate}
 \item The spinor $\phi$ has non-zero constant length and no zeros.
 \item The Ricci tensor decomposes as follows
\[{\rm Ric} 
\quad=\quad {1\over |\phi|^2}\sum_{k<l} \hat\eta_{kl}^\phi \circ \hat\Theta_{kl}
\quad=\quad {1\over |\phi|^2}\sum_{k<l} \hat\Theta_{kl} \circ \hat\eta_{kl}^\phi.
\]
 \item The scalar curvature is given by
\begin{eqnarray*}
{\rm R}
&=& {1\over |\phi|^2} \sum_{k<l} \tr(\hat\Theta_{kl} \circ\hat\eta_{kl}^\phi)
\quad=\quad -{1\over |\phi|^2}\left<\hat{\eta}^\phi,\hat{\Theta}\right>_1.
\end{eqnarray*}
 \item If the connection on the auxiliary bundle $F$ is flat, then $M$ is Ricci-flat. 
 \item If the parallel spinor $\phi$ is such that $\eta_{kl}^\phi=0$ for all $1\leq k< l\leq r$,
then
the manifold $M$ is Ricci-flat.
\end{enumerate}
\end{theo}

{\em Proof}.
Since the spinor $\phi$ is parallel
\begin{eqnarray*}
 X|\phi|^2 
 &=& \left< \nabla_X^\theta\phi, \phi\right>+\left< \phi, \nabla_X^\theta\phi\right>
 \quad=\quad 0.
\end{eqnarray*}
Thus, a non-trivial parallel spinor has constant length and no zeros.

Since $\phi$ is parallel, the left hand side of \rf{Ricci-curvature-identity}
is zero and
\begin{eqnarray*}
{\rm Ric}(e_j)\cdot\phi 
&=& \sum_{k<l}\sum_{s=1}^n\Theta_{kl}(e_j,e_s)e_s\cdot \kappa_{r*}^m(f_kf_l)\cdot\phi. 
\end{eqnarray*}
By taking the real part of the hermitian inner product with $e_i\cdot\phi$,
\begin{eqnarray*}
{\rm Re}\left<{\rm Ric}(e_j)\cdot\phi,e_i\cdot\phi\right> 
&=& \left<{\rm Ric}(e_j), e_i\right>|\phi|^2\\
&=& |\phi|^2{\rm Ric}_{ij}
\end{eqnarray*}
On the other hand,
\begin{eqnarray*}
\left<{\rm Ric}(e_j)\cdot\phi,e_i\cdot\phi\right> 
&=& \sum_{k<l}\sum_{s=1}^n\Theta_{kl}(e_j,e_s)\left<e_s\cdot
\kappa_{r*}^m(f_kf_l)\cdot\phi,e_i\cdot\phi\right>\\
&=& -\sum_{k<l}\sum_s\Theta_{kl}(e_j,e_s)\eta_{kl}^\phi(e_i,e_s)\\
&=& \sum_{k<l}(\hat\eta_{kl}^\phi\circ\hat\Theta_{kl})_{ij}.
\end{eqnarray*}
Hence, the Ricci endomorphism satisfies
\begin{eqnarray*}
|\phi|^2{\rm Ric} 
&=& \sum_{k<l} \hat\eta_{kl}^\phi \circ \hat\Theta_{kl} \\
&=& \sum_{k<l}   \hat\Theta_{kl} \circ \hat\eta_{kl}^\phi, 
\end{eqnarray*}
where the last equality is due to the symmetry of ${\rm Ric}$ and the skew-symmetry of both
$\hat\Theta_{kl}$ and $\hat\eta_{kl}^\phi$.
\qd

\subsubsection*{Example of a parallel twisted spinor}

Consider the subgroup
\[H:=\{[g,g]\in Spin(n)\times_{\mathbb{Z}_2}Spin(n)\,|\,\, g\in Spin(n)\}\subset
Spin(n)\times_{\mathbb{Z}_2}Spin(n).\] 
Clearly, $H$ is isomorphic to $SO(n)$, 
and the following diagram commutes
\[\begin{array}{ccc}
 &  & Spin(n)\times_{\mathbb{Z}_2}Spin(n)\\
 & \nearrow & \downarrow\\
H\cong SO(n) & \lra & SO(n)\times SO(n)
  \end{array}.
\]

\begin{prop}\label{universal-parallel-spinor}
 Every Riemannian manifold admits a spinorially twisted spin structure such that an associated
spinor bundle admits a parallel spinor field.
\end{prop}

{\em Proof}. Let $M$ be a Riemannian $n$-dimensional manifold.
Clearly, whether or not $M$ is Spin, it admits the twisted spin structure given by
\[\begin{array}{c}
P_{Spin^n(n)}\\
\downarrow\\
P_{SO(M)}\tilde{\times} P_{SO(M)}
  \end{array},
\]
where $P_{SO(M)}$ denotes the principal bundle of orthonormal frames of $M$.
Furthermore, by the diagram above we have a reduction of structure
\[\begin{array}{ccc}
&&P_{Spin^n(n)}\\
&\nearrow&\downarrow\\
P_{SO(M)}&\lra&P_{SO(M)}\tilde{\times} P_{SO(M)}
  \end{array},
\]
in such a way that the twisted spinor vector bundle $\Delta_n\otimes\Delta_n$ can be seen as an
associated vector bundle to $P_{SO(M)}$, which is a well known fact.

Let $\mathcal{B}$ be the unitary basis of $\Delta_n$ described in Section \ref{preliminaries} and
$\gamma_n$ be the corresponding real or quaternionic structure of $\Delta_n$. We claim that the spinor
\begin{eqnarray*}
\phi_0
&:=&\sum_{\psi\in\mathcal{B}} \psi\otimes \gamma_n(\psi)\\
&=&\sum_{(\varepsilon_1,\ldots,\varepsilon_{[n/2]})\in\{1,-1\}^{\times[n/2]}}
C(n,\varepsilon_1,\ldots,\varepsilon_{[n/2]})\,\,
u_{\varepsilon_1,\ldots,\varepsilon_{[n/2]}}\otimes
u_{-\varepsilon_1,\ldots,-\varepsilon_{[n/2]}} 
\end{eqnarray*}
is invariant under $H\cong SO(n)$, where 
\[
\begin{array}{ccll}
 C(n;\varepsilon_1,\ldots,\varepsilon_{4k}) 
&=&
(-1)^{k+{1\over 2}\sum_{j=1}^{2k} (\varepsilon_{2j-1}+1)}
 &\mbox{if $n=8k,8k+1$,} \\
 C(n;\varepsilon_1,\ldots,\varepsilon_{4k+1}) 
&=&
i(-1)^{k+{1\over 2}\sum_{j=1}^{2k+1} (\varepsilon_{2j-1}+1)}
 &\mbox{if $n=8k+2,8k+3$,} \\
 C(n;\varepsilon_1,\ldots,\varepsilon_{4k+2}) 
&=&
(-1)^{k+{1\over 2}\sum_{j=1}^{2k+1} (\varepsilon_{2j-1}+1)}
 &\mbox{if $n=8k+4,8k+5$,} \\
 C(n;\varepsilon_1,\ldots,\varepsilon_{4k+3}) 
&=&
i(-1)^{k+{1\over 2}\sum_{j=1}^{2k+2} (\varepsilon_{2j-1}+1)}
 &\mbox{if $n=8k+6,8k+7$.} 
\end{array}
\]
We will prove the invariance by means of the Lie algebra
$Lie(H)\cong\mathfrak{so}(n)$. 
Let us consider the case
$n=8k$.
Let $\{e_1,\ldots, e_{8k}\}\subset\mathbb{R}^{8k}$ be an ordered orthonormal basis, and
$\{f_1,\ldots, f_{8k}\}\subset\mathbb{R}^{8k}$ be the same ordered basis but renamed since it will
refer to the auxiliary twisting bundle. Thus
\[Lie(H) = \span\{e_ie_j + f_if_j \in\mathfrak{spin}(n)\oplus\mathfrak{spin}(n)\, |\,\, 1\leq
i<j\leq n\}.\]
Let us consider one summand in $\phi_0$,
\[\varphi_1:=(-1)^{{1\over 2}\sum_{j=1}^{2k}
\varepsilon_{2j-1}}\,\,u_{\varepsilon_1,\ldots,\varepsilon_{4k}}\otimes
u_{-\varepsilon_1,\ldots,-\varepsilon_{4k}},\]
and focus first on 
\[u_{\varepsilon_1,\ldots,\varepsilon_{4k}}\otimes
u_{-\varepsilon_1,\ldots,-\varepsilon_{4k}},\]
Recall that for $1\leq p\leq 8k$,
\begin{eqnarray*}
e_{p}\cdot u_{\varepsilon_1,\ldots,\varepsilon_{4k}}
&=& i^{p-2[p/2]}(-1)^{[(p+1)/2]-1}
\left(\prod_{\alpha=4k-[(p+1)/2]+1+p-2[p/2]}^{4k} \varepsilon_{\alpha}\right) 
u_{\varepsilon_1,\ldots, (-\varepsilon_{4k-[(p+1)/2]+1}) ,\ldots,\varepsilon_{4k}}.    
\end{eqnarray*}
If we apply $e_pe_q$ to it, with $1\leq p< q\leq n$ and $[(q-1)/2]>[(p-1)/2]$, we get
\begin{eqnarray*}
&&i^{q-2[q/2]}(-1)^{[(q+1)/2]-1}
\left(\prod_{\alpha=4k-[(q+1)/2]+1+q-2[q/2]}^{4k} \varepsilon_{\alpha}\right) \\
&&\times  
i^{p-2[p/2]}(-1)^{[(p+1)/2]-1}
\left(\prod_{\beta=4k-[(p+1)/2]+1+p-2[p/2]}^{4k} \varepsilon_{\beta}\right) \\
&&\times  u_{\varepsilon_1,\ldots, (-\varepsilon_{4k-[(q+1)/2]+1}) 
,\ldots,(-\varepsilon_{4k-[(p+1)/2]+1}),\ldots,\varepsilon_{4k}}
\otimes
u_{-\varepsilon_1,\ldots,-\varepsilon_{4k}}. 
\end{eqnarray*}
Now, let us consider another summand in $\phi_0$ 
\begin{eqnarray*}
\varphi_2 &:=& C(8k,\varepsilon_1,\ldots, (-\varepsilon_{4k-[(q+1)/2]+1}) 
,\ldots,(-\varepsilon_{4k-[(p+1)/2]+1}),\ldots,\varepsilon_{4k}) \\
&&\times\,\,u_{\varepsilon_1,\ldots, (-\varepsilon_{4k-[(q+1)/2]+1}) 
,\ldots,(-\varepsilon_{4k-[(p+1)/2]+1}),\ldots,\varepsilon_{4k}}
\otimes
u_{-\varepsilon_1,\ldots, (\varepsilon_{4k-[(q+1)/2]+1}) 
,\ldots,(\varepsilon_{4k-[(p+1)/2]+1}),\ldots,-\varepsilon_{4k}},
\end{eqnarray*}
and focus first on 
\[u_{\varepsilon_1,\ldots, (-\varepsilon_{4k-[(q+1)/2]+1}) 
,\ldots,(-\varepsilon_{4k-[(p+1)/2]+1}),\ldots,\varepsilon_{4k}}
\otimes
u_{-\varepsilon_1,\ldots, (\varepsilon_{4k-[(q+1)/2]+1}) 
,\ldots,(\varepsilon_{4k-[(p+1)/2]+1}),\ldots,-\varepsilon_{4k}}.\]
When we apply $\kappa_{n*}^1(f_pf_q)$ to it,
\begin{eqnarray*}
 &&  
i^{q-2[q/2]}(-1)^{[(q+1)/2]-1}
\left(\prod_{\alpha=4k-[(q+1)/2]+1+q-2[q/2]}^{4k} -\varepsilon_{\alpha}\right)  \\ 
&&\times 
i^{p-2[p/2]}(-1)^{[(p+1)/2]-1}
\left(\prod_{\beta=4k-[(p+1)/2]+1+p-2[p/2]}^{4k} -\varepsilon_{\beta}\right) \\
&&\times (-1)(-1)^{q+1}(-1)^{p+1}
\\
&&\times  u_{\varepsilon_1,\ldots, (-\varepsilon_{4k-[(q+1)/2]+1}) 
,\ldots,(-\varepsilon_{4k-[(p+1)/2]+1}),\ldots,\varepsilon_{4k}}
\otimes
u_{-\varepsilon_1,\ldots,-\varepsilon_{4k}}.
\end{eqnarray*}
Now, while the coefficient of $e_pe_q\cdot \varphi_1$ is
\begin{eqnarray*}
&& C(8k;\varepsilon_1,\ldots,\varepsilon_{4k})\\
&&\times\,\,i^{q-2[q/2]}(-1)^{[(q+1)/2]-1}
\left(\prod_{\alpha=4k-[(q+1)/2]+1+q-2[q/2]}^{4k} \varepsilon_{\alpha}\right) \\
&&\times  
i^{p-2[p/2]}(-1)^{[(p+1)/2]-1}
\left(\prod_{\beta=4k-[(p+1)/2]+1+p-2[p/2]}^{4k} \varepsilon_{\beta}\right),
\end{eqnarray*}
the coefficient of $\kappa_{n*}^1(f_pf_q)\cdot \varphi_2$ is
\begin{eqnarray*}
&&C(8k,\varepsilon_1,\ldots, (-\varepsilon_{4k-[(q+1)/2]+1}) 
,\ldots,(-\varepsilon_{4k-[(p+1)/2]+1}),\ldots,\varepsilon_{4k}) \\
&&\times (-1)^{1+[(q+1)/2]+[(p+1)/2]}\\
&&\times\,\,
i^{q-2[q/2]}(-1)^{[(q+1)/2]-1}
\left(\prod_{\alpha=4k-[(q+1)/2]+1+q-2[q/2]}^{4k} \varepsilon_{\alpha}\right)  \\ 
&&\times
i^{p-2[p/2]}(-1)^{[(p+1)/2]-1}
\left(\prod_{\beta=4k-[(p+1)/2]+1+p-2[p/2]}^{4k} \varepsilon_{\beta}\right).
\end{eqnarray*}
By checking the possible cases in which $[(p+1)/2]$ and $[(q+1)/2]$ are either even or odd,
these two coefficients differ by $(-1)$. 
Thus
\[e_pe_q\cdot\varphi_1 + \kappa_{n*}^1(f_pf_q)\cdot\varphi_2 =0.\]
Clearly, every summand in $\phi_0$ has a unique counterpart as in the previous calculation.
All the other possible cases for values and parities of $n$, $p$ and $q$ are similar.
Hence $Lie(H)\cong\mathfrak{so}(n)$ annihilates $\phi_0$.
\qd

\begin{prop}\label{second-outstanding-feature}
The 2-forms associated to $\phi_0$ are 
multiples of the elements of the stanfard basis of $\mathfrak{so}(n)$, i.e.
\[\eta_{pq}^{\phi_0} = 2^{[n/2]}\,\,  e_p\wedge e_q.\]
\end{prop}

{\em Proof}.
Notice that
\begin{eqnarray*}
\phi_0
&=&\sum_{(\varepsilon_1,\ldots,\varepsilon_{[n/2]})\in\{1,-1\}^{\times[n/2]}}
C(n,\varepsilon_1,\ldots,\varepsilon_{[n/2]})\,\,
u_{\varepsilon_1,\ldots,\varepsilon_{[n/2]}}\otimes
u_{-\varepsilon_1,\ldots,-\varepsilon_{[n/2]}}. 
\end{eqnarray*}
is orthogonal to any spinor orthogonal to
\[V_0=\span\{u_{\varepsilon_1,\ldots,\varepsilon_{[n/2]}}\otimes
u_{-\varepsilon_1,\ldots,-\varepsilon_{[n/2]}} \, | \,\, 
(\varepsilon_1,\ldots,\varepsilon_{[n/2]})\in\{1,-1\}^{\times[n/2]} \}.\]
Thus, for $p<q$, $s<t$, $(p,q)\not = (s,t)$,
\begin{eqnarray*}
\eta_{st}^{\phi_0}(e_p,e_q)&=&\left<e_pe_q\cdot \kappa_{n*}^1(f_sf_t)\cdot \phi_0 , \phi_0\right> \\
&=& 0,
\end{eqnarray*}
since each one of the summands in $e_pe_q\cdot \kappa_{n*}^1(f_sf_t)\cdot \phi_0$ is orthogonal to
$V_0$.

On the other hand, if $(p,q) = (s,t)$ with $1\leq p< q\leq n$, $[(q-1)/2]>[(p-1)/2]$, and
\[\varphi_1:=(-1)^{{1\over 2}\sum_{j=1}^{2k}
\varepsilon_{2j-1}}\,\,u_{\varepsilon_1,\ldots,\varepsilon_{4k}}\otimes
u_{-\varepsilon_1,\ldots,-\varepsilon_{4k}},\]
then $ e_pe_q\cdot \kappa_{n*}^1(f_pf_q)\cdot \varphi_1 $ is equal to
\begin{eqnarray*}
&&(-1)^{[(p+1)/2]+[(q+1)/2]+{1\over 2}\sum_{j=1}^{2k}\varepsilon_{2j-1}}\\
&&\times  u_{\varepsilon_1,\ldots, (-\varepsilon_{4k-[(q+1)/2]+1}) 
,\ldots,(-\varepsilon_{4k-[(p+1)/2]+1}),\ldots,\varepsilon_{4k}}
\otimes
u_{-\varepsilon_1,\ldots, (\varepsilon_{4k-[(q+1)/2]+1}) 
,\ldots,(\varepsilon_{4k-[(p+1)/2]+1}),\ldots,-\varepsilon_{4k}}.
\end{eqnarray*}
Clearly, this is paired with
\begin{eqnarray*}
\varphi_3&=& C(8k;\varepsilon_1,\ldots, (-\varepsilon_{4k-[(q+1)/2]+1}) 
,\ldots,(-\varepsilon_{4k-[(p+1)/2]+1}),\ldots,\varepsilon_{4k})\\
&& u_{\varepsilon_1,\ldots, (-\varepsilon_{4k-[(q+1)/2]+1}) 
,\ldots,(-\varepsilon_{4k-[(p+1)/2]+1}),\ldots,\varepsilon_{4k}}
\otimes
u_{-\varepsilon_1,\ldots, (\varepsilon_{4k-[(q+1)/2]+1}) 
,\ldots,(\varepsilon_{4k-[(p+1)/2]+1}),\ldots,-\varepsilon_{4k}}
\end{eqnarray*}
in the hermitian product, so that 
\begin{eqnarray*}
\left<e_pe_q\cdot \kappa_{n*}^1(f_pf_q)\cdot\varphi_1, \varphi_3\right>
&=&1,
\end{eqnarray*}
for all possible cases in which $[(p+1)/2]$ and $[(q+1)/2]$ are either even or odd.
Furthermore, all the other cases for values and parities of $n$, $p$ and $q$ are
similar.
Since $\phi_0$ is made up of $2^{[n/2]}$ summands which satisfy the previous arguments, 
\[\eta_{pq}^{\psi_0}(e_s,e_t)=2^{[n/2]}\,(\delta_{ps}\delta_{qt}-\delta_{pt}\delta_{qs}).\]
\qd

Let us now check our curvature formulas on this example.
Formula \rf{curvature-twisted-spin} becomes
\begin{eqnarray*}
R^{\theta}(X,Y)(\phi_0)&=&
{1\over 2} \sum_{1\leq i<j\leq n} \Omega_{ij}(X,Y) e_ie_j\cdot\phi_0+
{1\over 2} \sum_{1\leq k<l\leq n} \Theta_{kl}(X,Y) \kappa_{n*}^1(f_kf_l)\cdot\phi_0\\
&=&
{1\over 2} \sum_{i<j} \big<R^M(X,Y)(e_i),e_j \big> e_ie_j\cdot\phi_0+
{1\over 2} \sum_{i<j} \big<R^M(X,Y)(e_i),e_j \big>
\kappa_{n*}^1(f_if_j)\cdot\phi_0\\
&=&
{1\over 2} \sum_{i<j} \big<R^M(X,Y)(e_i),e_j \big> (e_ie_j+\kappa_{n*}^1(f_if_j))\cdot\phi_0\\
&=&0,
\end{eqnarray*}
which is consistent with the parallelness of $\phi_0$, 
and
\begin{eqnarray*}
\sum_{1\leq k<l\leq n} (\hat\Omega_{kl} \circ\hat\eta_{kl}^{\phi_0})_{st}
&=& \sum_{1\leq k<l\leq n}\sum_{a=1}^n (\hat\Omega_{kl})_{sa} (\hat\eta_{kl}^{\phi_0})_{at}\\ 
&=& \sum_{a=1}^n\sum_{1\leq k<l\leq n} \Omega_{kl}(e_s,e_a) \eta_{kl}^{\phi_0}(e_a,e_t)\\ 
&=& 2^{[n/2]}\sum_{a=1}^n\sum_{1\leq k<l\leq n} \left<R(e_s,e_a)e_k,e_l\right>
(\delta_{ka}\delta_{lt}-\delta_{kt}\delta_{la})\\ 
&=& 2^{[n/2]}\left(\sum_{a<t}
\left<R(e_s,e_a)e_a,e_t\right>
-\sum_{a>t}\left<R(e_s,e_a)e_t,e_a\right>\right)\\ 
&=& 2^{[n/2]}\sum_{a}
\left<R(e_a,e_s)e_t,e_a\right>\\ 
&=& 2^{[n/2]}\sum_{a}
{\rm Ric}_{st}.
\end{eqnarray*}
which is consistent with \rf{Ricci-curvature-identity}.

\subsection{Killing spinors}\label{section-killing-spinors}

\begin{defi}
A spinor $\phi\in\Gamma(S)$ is a Killing spinor if, for every $X\in \Gamma(TM)$,
\[\nabla_X\phi = \mu\, X\cdot\phi,\]
where $\mu\in\mathbb{C}$.
\end{defi}

\begin{prop} Let $\phi\in\Gamma(S)$ be a non-trivial Killing spinor.
\begin{enumerate}
 \item $\phi$ has no zeros.
 \item $\phi$ is an eigenspinor of the twisted Dirac operator.
 \item If $\phi$ is a {\em real Killing spinor}, i.e. $\mu$ is real, the length of the Killing
spinor $\phi$ is constant.
 \item If $\mu$ is real, then the vector field
\[V^\phi = \sum_n \left<e_i\cdot \phi,\phi\right>e_i\]
is a Killing vector field.
\end{enumerate}
\end{prop}
The proof is analogous to the one in the Spin$^c$ case (cf. \cite{friedrich}). \qd

\begin{theo}\label{Ricci-Killing-spinor}
 Let $\phi\in\Gamma(S)$ be a real Killing spinor.
Then 
\begin{itemize}
 \item The Ricci tensor decomposes as follows
\[{\rm Ric} = 4(n-1)\mu^2\, {\rm Id}_{TM} + {1\over |\phi|^2}\sum_{k<l} \hat\Theta_{kl}
\circ\hat\eta_{kl}^\phi.\]
 \item The scalar curvature is given by
\begin{eqnarray*}
{\rm R}
&=&  4n(n-1)\mu^2+{1\over |\phi|^2}\sum_{k<l} \tr(\hat\Theta_{kl} \circ\hat\eta_{kl}^\phi). 
\end{eqnarray*}
 \item If the connection on the auxiliary bundle $F$ is flat, then $M$ is Einstein. 
 \item If the Killing spinor $\phi$ is such that $\eta_{kl}^\phi=0$ for all $1\leq k< l\leq r$,
then
the manifold $M$ is Einstein.
\end{itemize}
\end{theo}
{\em Proof}.
The left hand side of \rf{Ricci-curvature-identity} now becomes
\begin{eqnarray*}
 \sum_{i = 1}^n e_i\cdot R^{\theta}(e_j, e_i)(\phi) 
&=& 
 \sum_{i \not= j} e_i\cdot \mu^2 (e_ie_j-e_je_i)\cdot\phi\\
&=& 
 -2\mu^2\sum_{i \not= j}  e_j\cdot\phi\\
&=& 
 -2(n-1)\mu^2  e_j\cdot\phi.
\end{eqnarray*}
By taking the real part of the hermitian product with $e_t\cdot\phi$
we get
\begin{eqnarray*}
 {\rm Re}\left[-2(n-1)\mu^2 \left<e_j\cdot \phi,e_t\cdot \phi\right>\right] 
&=& -2(n-1)\mu^2 \left<e_j,e_t\right>|\phi|^2\\
&=& -2|\phi|^2(n-1)\mu^2 \delta_{jt}.
\end{eqnarray*}
Hence, by the calculations of the last subsection we have
\[{\rm Ric} = 4(n-1)\mu^2\, {\rm Id}_{TM} + {1\over |\phi|^2}\sum_{k<l} \hat\Theta_{kl}
\circ\hat\eta_{kl}^\phi.\]
\qd

\subsection{Generalized real Killing spinors}

\begin{defi}
A spinor $\phi\in\Gamma(S)$ is called a {\em generalized Killing spinor} if
\[\nabla_X\phi=-E(X)\cdot\phi\]
for some symmetric endomorphism $E\in\Gamma(\End(TM))$ and all $X\in\Gamma(TM)$. 
\end{defi}

In this case, the left hand side of \rf{Ricci-curvature-identity} is 
\begin{eqnarray*}
 \sum_{i = 1}^n e_i\cdot R^{\theta}(e_s, e_i)(\phi) 
&=&
 \sum_{i = 1}^n e_i\cdot R^{\theta}(e_s, e_i)(\phi) \\
&=&
 \sum_{i =1}^n e_i\cdot
(\nabla^\theta_{e_s}\nabla^\theta_{e_i}-\nabla^\theta_{e_i}\nabla^\theta_{e_s}-\nabla^\theta_{[e_s,
e_i]} )\phi \\
&=&
 \sum_{i =1}^n e_i\cdot 
((\nabla_{e_i}E)(e_s)-(\nabla_{e_s}E)(e_i)+E(e_i)\cdot E(e_s)-E(e_s)\cdot
E(e_i))\cdot\phi    \\
&=&
 \sum_{i \not=s} e_i\cdot 
(d^\nabla E(e_i,e_s)+E(e_i)\cdot E(e_s)-E(e_s)\cdot E(e_i))\cdot\phi ,   
\end{eqnarray*}
where
\[d^\nabla E(X,Y)=(\nabla_{X}E)(Y)-(\nabla_{Y}E)(X).\]
Now, if the orthonormal frame also diagonalizes $E$, for $i\not=s$,
\begin{eqnarray*}
 E(e_i)\cdot E(e_s)-E(e_s)\cdot E(e_i)
&=& E_{ii}E_{ss}(e_i\cdot e_s - e_s\cdot e_i)\\
&=& 2E_{ii}E_{ss}e_i\cdot e_s,
\end{eqnarray*}
and
\begin{eqnarray*}
e_i\cdot( E(e_i)\cdot E(e_s)-E(e_s)\cdot E(e_i))
&=& -2E_{ii}E_{ss} e_s, 
\end{eqnarray*}
so that
\begin{eqnarray*}
\sum_{i\not=s} e_i\cdot( E(e_i)\cdot E(e_s)-E(e_s)\cdot E(e_i))
&=& -2\left(\sum_{i\not=s}E_{ii} \right)E_{ss}e_s\\ 
&=& -2\tr(E) E_{ss}e_s +2 E_{ss}^2 e_s\\ 
&=& 2E_{ss}(E_{ss}-\tr(E))  e_s. 
\end{eqnarray*}
By taking the real part of the hermitian product with $e_t\cdot\phi$
\begin{eqnarray*}
{\rm Re}\left<\sum_{i\not=s} e_i\cdot( E(e_i)\cdot E(e_s)-E(e_s)\cdot E(e_i))\cdot
\phi,e_t\cdot\phi\right>
&=& {\rm Re}\left<2E_{ss}(E_{ss}-\tr(E))  e_s\cdot\phi,e_t\cdot\phi\right>\\ 
&=& 2E_{ss}(E_{ss}-\tr(E)){\rm Re}\left<  e_s\cdot\phi,e_t\cdot\phi\right>\\ 
&=& 2E_{ss}(E_{ss}-\tr(E))\delta_{ts} |\phi|^2,
\end{eqnarray*}
which gives the matrix
\[2|\phi|^2(E^2-\tr(E)E).\]

On the other hand,
\begin{eqnarray*}
  \sum_{i \not=s} e_i\cdot d^\nabla E(e_i,e_s) 
&=&  \sum_{i \not=s} (e_i\wedge d^\nabla E(e_i,e_s) - \left<e_i, d^\nabla E(e_i,e_s)\right> )\\
&=&  \sum_{i \not=s} \left(e_i\wedge \left(\sum_{j=1}^n\left<d^\nabla E(e_i,e_s),e_j\right>
e_j\right) -
\left<e_i, d^\nabla E(e_i,e_s)\right> \right)\\
&=&  \sum_{i \not=s} \left(\sum_{j=1}^n\left<d^\nabla E(e_i,e_s),e_j\right>e_i\wedge  e_j -
\left<e_i, d^\nabla E(e_i,e_s)\right> \right)\\
&=&  \sum_{i \not=s} \left(\sum_{j\not=i}\left<d^\nabla E(e_i,e_s),e_j\right>e_i\cdot  e_j -
\left<e_i, d^\nabla E(e_i,e_s)\right> \right)
\end{eqnarray*}
By taking the real part of the hermitian product with $e_t\cdot\phi$
we get
\begin{eqnarray}
 {\rm Re}\left< \sum_{i \not=s} e_i\cdot d^\nabla E(e_i,e_s) \cdot \phi,e_t\cdot\phi\right>
&=&  {\rm Re}\left<\sum_{i \not=s} \left(\sum_{j\not=i}\left<d^\nabla E(e_i,e_s),e_j\right>e_i\cdot 
e_j -\left<e_i, d^\nabla E(e_i,e_s)\right> \right)\cdot\phi,e_t\cdot\phi\right>\nonumber\\
&=&  {\rm Re}\left<\sum_{i \not=s} \left(\sum_{j\not=i}\left<d^\nabla E(e_i,e_s),e_j\right>e_i\cdot 
e_j\right)\cdot\phi,e_t\cdot\phi\right>\nonumber\\
&=& - \sum_{i \not=s}\sum_{j\not=i}\left<d^\nabla E(e_i,e_s),e_j\right>{\rm Re}\left< e_t\cdot
e_i\cdot e_j \cdot\phi,\phi\right>\nonumber\\
&=&  \sum_{i,j}\left<d^\nabla E(e_s,e_i),e_j\right>\nu^\phi(e_t, e_i,e_j) \nonumber\\
&=&  \sum_{i,j}\left<(e_s\lrcorner d^\nabla E)(e_i),e_j\right>(e_t\lrcorner\nu^\phi)(e_i,e_j) \nonumber\\
&=&  \sum_{i,j}(e_s\lrcorner d^\nabla E)_{ji}(e_t\lrcorner\nu^\phi)_{ji} \nonumber\\
&=:&  ((\lrcorner d^\nabla E)\circledast(\lrcorner\nu^\phi))_{ts},\label{eq:circulo-asterisco}
\end{eqnarray}
where
\begin{eqnarray}
 \nu^\phi(e_t, e_i,e_j) 
&:=& {\rm Re}\left< e_t\cdot e_i\cdot e_j \cdot\phi,\phi\right>  \label{eq:def-nu}. 
\end{eqnarray}
Thus,
\begin{eqnarray*}
 {\rm Ric}_{st}
&=& 
-2E_{ss}(E_{ss}-\tr(E))\delta_{st} - \left({\lrcorner
d^\nabla E}\circledast{\lrcorner\nu^\phi}\right)_{st}+
\sum_{k<l}(\hat\Theta_{kl}\circ\hat\eta_{kl}^\phi)_{st}, 
\end{eqnarray*}
i.e.
\[{\rm Ric} = -2E^2+2\tr(E)E - \left({\lrcorner
d^\nabla E}\circledast{\lrcorner\nu^\phi}\right)  + \sum_{k<l}\hat\Theta_{kl}\circ\hat\eta_{kl}^\phi.\]
and
\[{\rm R} =  -2\tr(E^2)+2\tr(E)^2- \tr\left({\lrcorner
d^\nabla E}\circledast{\lrcorner\nu^\phi}\right) 
+\sum_{k<l}\tr(\hat\Theta_{kl}\circ\hat\eta_{kl}^\phi).\]
Thus, we have proved the following.

\begin{theo}
 Let $\phi\in\Gamma(S)$ be a generalized Killing spinor.
Then 
\begin{itemize}
 \item the Ricci tensor decomposes as follows
\[{\rm Ric} = -2E^2+2\tr(E)E - \left({\lrcorner
d^\nabla E}\circledast{\lrcorner\nu^\phi}\right)  + \sum_{k<l}\hat\Theta_{kl}\circ\hat\eta_{kl}^\phi,\]
where  $\left({\lrcorner d^\nabla E}\circledast{\lrcorner\nu^\phi}\right)$ and $\nu^\phi$ are defined as in
{\em \rf{eq:circulo-asterisco}} and {\em \rf{eq:def-nu}} respectively;
 \item the scalar curvature is given by
\[{\rm R} =  -2\tr(E^2)+2\tr(E)^2- \tr\left({\lrcorner
d^\nabla E}\circledast{\lrcorner\nu^\phi}\right) 
+\sum_{k<l}\tr(\hat\Theta_{kl}\circ\hat\eta_{kl}^\phi).\]
\end{itemize}
\end{theo}
\qd

{\bf Remark}.
These formulas reduce to the previous two cases when
$E$ is a multiple of the identity endomorphism
$E=\mu\,\, {\rm Id}_{TM}$.

\subsection{Twisted Schr\"odinger-Lichnerowicz formula}\label{section-twisted-SL-formula}

Recall the curvature operator 
\[\Theta = \sum_{k<l} \Theta_{kl} \otimes f_kf_l \in \ext^2 TM \otimes \ext^2 F\]
of the connection $F$, and denote by
\begin{eqnarray*}
\tilde\Theta^m
&=& (\mu_n\otimes \kappa_{r*}^m) (\Theta) 
\end{eqnarray*}
the corresponding operator on twisted spinor fields. 
For future use, note the following operator norm inequality
\[|\tilde\Theta^m| \leq m|\tilde\Theta^1|,\]
which follows from \rf{kappa-en-mus}.

\begin{theo}[Twisted Schr\"odinger-Lichnerowicz Formula]\label{theo-SL-formula}
Let $\phi\in\Gamma(S)$. Then
\begin{equation}
 \dirac^\theta(\dirac^\theta(\phi)) =  \Delta^{\theta}(\phi)+ \dfrac{\rm R}{4}\phi +
\dfrac{1}{2}\tilde\Theta^m\cdot\phi \label{SL-formula}
\end{equation}
where ${\rm R}$ is the scalar curvature of the Riemannian manifold $M$.
\end{theo}
{\em Proof}. 
Consider the
difference
\begin{eqnarray*}
\dirac^\theta(\dirac^\theta(\phi))  - \Delta^{\theta}(\phi)
&=& \sum_{i,j} e_i\cdot\nabla_{e_i}^{\theta}(e_j\cdot\nabla_{e_j}^{\theta}\phi)
+\sum_i \nabla_{e_i}^{\theta}\nabla_{e_i}^{\theta}\phi + \sum_i
\diver(e_i)\nabla_{e_i}^{\theta}\phi\\
&=& \sum_{i,j} e_i\cdot(\nabla_{e_i}e_j\cdot\nabla_{e_i}^{\theta}\phi
+e_j\cdot\nabla_{e_j}^{\theta}\nabla_{e_j}^{\theta}\phi)
+\sum_i \nabla_{e_i}^{\theta}\nabla_{e_i}^{\theta}\phi + \sum_i
\diver(e_i)\nabla_{e_i}^{\theta}\phi\\
&=& \sum_{i,j,k} \left<\nabla_{e_i}e_j,e_k\right>e_ie_k\cdot\nabla_{e_i}^{\theta}\phi
+\sum_{i,j}e_ie_j\cdot\nabla_{e_j}^{\theta}\nabla_{e_j}^{\theta}\phi)\\
&&+\sum_i \nabla_{e_i}^{\theta}\nabla_{e_i}^{\theta}\phi + \sum_i
\diver(e_i)\nabla_{e_i}^{\theta}\phi\\
&=& \sum_{i}\sum_{j\not =k} \left<\nabla_{e_i}e_j,e_k\right>e_ie_k\cdot\nabla_{e_i}^{\theta}\phi
+\sum_{i\not=j}e_ie_j\cdot\nabla_{e_j}^{\theta}\nabla_{e_j}^{\theta}\phi),
\end{eqnarray*}
since
\[\sum_j\sum_{i=k}\left<\nabla_{e_i}e_j,e_k\right>e_ie_k \nabla_{e_j}^{\theta}\phi =
-\sum_j\diver(e_j)\nabla_{e_j}^{\theta}\phi.\]
Now, for fixed $j$
\[\sum_{i\not=k}\left<\nabla_{e_i}e_j,e_k\right>e_ie_k
=\sum_{i<k}\left<e_j,[e_k,e_i]\right>e_ie_k.\]
Thus,
\begin{eqnarray*}
\dirac^\theta(\dirac^\theta(\phi))  - \Delta^{\theta}(\phi)
&=& \sum_{j}\sum_{i<k} \left<e_j,[e_k,e_i]\right>e_ie_k\cdot\nabla_{e_i}^{\theta}\phi
+\sum_{i<j}e_ie_j\cdot(\nabla_{e_i}^{\theta}\nabla_{e_j}^{\theta}
-\nabla_{e_j}^{\theta}\nabla_{e_i}^{ \theta})\phi\\
&=& \sum_{i<j}e_ie_j\cdot(\nabla_{e_i}^{\theta}\nabla_{e_j}^{\theta}
-\nabla_{e_j}^{\theta}\nabla_{e_i}^{ \theta}-\nabla_{[e_i,e_j]}^{\theta})\phi\\
&=&{1\over 2}\sum_{i,j}e_ie_jR^\theta(e_i,e_j)\phi.
\end{eqnarray*}
The result follows from Proposition \ref{scalar-curvature-identity}.
\qd

\subsection{Bochner-type arguments}

In this subsection we will prove some corollaries of the Schr\"odinger-Lichnerowicz formula and
Bochner type arguments as in \cite{friedrich,Hijazi, Nakad}.
From here onwards, we will assume that the $n$-dimensional Riemannian Spin$^r$ manifold $M$ is compact
(without border) and connected.

\subsubsection{Harmonic spinors}

\begin{corol}\label{corol-no-harmonic-spinors} 
If ${\rm R}\geq 2m|\tilde\Theta^1|$ everywhere (in point-wise operator norm), and the inequality is
strict at a point, then
\[\ker(\dirac^{\theta,m}) =0.\]
Furthermore,
\[\ker(\dirac^{\theta,m'}) =0\]
for any $0\leq m'\leq m$ such that the bundle $S(TM)\otimes S(F)^{\otimes m'}$ is globally defined.
\end{corol}

{\em Proof}. If $\phi\not=0$ is a solution of 
\[\dirac^\theta(\phi) =0,\]
by the twisted Schr\"odinger-Lichnerowicz formula \rf{SL-formula}
\[0 =   \Delta_{\theta}(\phi)+ \dfrac{\rm R}{4}\phi +
\dfrac{1}{2}\tilde\Theta^m\cdot\phi.\]
By taking  hermitian product with $\phi$ and integrating over $M$ we get
\begin{eqnarray*}
0
&=&  
\int_M\left< \Delta^{\theta}(\phi),\phi\right>+ \int_M\dfrac{\rm R}{4}\left<\phi,\phi\right> +
\dfrac{1}{2}\int_M\left<\tilde\Theta^m\cdot\phi,\phi\right>\\ 
&\geq&  
\int_M| \nabla^\theta\phi|^2+ \dfrac{1}{4}\int_M\left({\rm R} -
2|\tilde\Theta^m|\right)|\phi|^2\\
&\geq&  
\int_M| \nabla^\theta\phi|^2+ \dfrac{1}{4}\int_M\left({\rm R} -
2m|\tilde\Theta^1|\right)|\phi|^2.
\end{eqnarray*}
Since
\[{\rm R} -
2m|\tilde\Theta^1|\geq 0,\]
then 
\[|\nabla^\theta\phi|=0,\]
so that $\phi$ is parallel,
has non-zero constant length and 
no zeros.
Furthermore, since 
\[{\rm R} -
2m|\tilde\Theta^1| > 0\]
at some point,
\[0
\geq  
 |\phi|^2\int_M\left({\rm R} -
2m|\tilde\Theta^1|\right)>0,
\]
which is a contradiction.

The last claim now follows from 
\[{\rm R}\geq {\rm R}- 2|\tilde\Theta^1|\geq {\rm R}-4|\tilde\Theta^1|\geq\cdots\geq {\rm
R}-2m|\tilde\Theta^1|\geq 0.\]
\qd

{\bf Remarks}.
\begin{itemize}
 \item The last statement of Corollary \ref{corol-no-harmonic-spinors} means that if there are no harmonic
spinors for a given power due to the condition on the curvatures ${\rm R}$ and $\Theta^1$, then there are no
harmonic spinors for the twisted spinor bundles with smaller powers of $\Delta_r$ either.
 \item One can prove that a compact Riemannian $n$-dimensional manifold carrying a non-flat
parallel even Clifford structure of rank $r$ (cf. \cite{Moroianu-Semmelmann}) is a Spin$^r$ manifold and
carries no harmonic spinors for 
\begin{equation}
m\leq {n+8r-16\over r(r-1)}\label{eq:power-upper-bound} 
\end{equation}
if the scalar curvature is non-negative.

In particular, the case of rank $r=3$ corresponds to quaternion-K\"ahler manifolds, and Corollary
\ref{corol-no-harmonic-spinors}
and \rf{eq:power-upper-bound}
reproduce some of the vanishings of indices of twisted Dirac operators proved (via twistor transform) in
\cite{Salamon}.

\end{itemize}

\vspace{.1in}

Now notice that
\begin{eqnarray}
\left<\tilde\Theta^m\cdot\phi,\phi\right> 
&=& \left<\sum_{k<l} \left[\sum_{i<j}\Theta_{kl}(e_i,e_j)e_ie_j\right]\cdot\kappa_{r*}^m(f_kf_l)
\cdot\phi,\phi\right>\nonumber\\
&=& \sum_{k<l} \sum_{i<j}\Theta_{kl}(e_i,e_j)\eta_{kl}^\phi(e_i,e_j)\nonumber\\
&=& \left<{\Theta},{\eta}^{\phi}\right>_0,\label{real-number}
\end{eqnarray}
which is a real number dependent on the curvature of the connection on $F$ and on the specific
spinor $\phi$.

\begin{corol}\label{not-massless-Dirac-spinor}
If $\phi$ is such that
\[{\rm R}|\phi|^2+ 2\left<{\Theta},{\eta}^\phi\right>_0\geq 0\] 
everywhere, and the inequality
is strict at a
point, then 
\[\dirac^\theta(\phi) \not=0.\]
\end{corol}
{\em Proof}. Suppose $\phi\not=0$ is such that 
\[\dirac^\theta(\phi) =0.\]
Then, by \rf{SL-formula}
\begin{eqnarray*}
0
&=&  
\int_M\left< \Delta^{\theta}(\phi),\phi\right>+ \int_M\dfrac{\rm R}{4}\left<\phi,\phi\right> +
\dfrac{1}{2}\int_M\left<\tilde\Theta\cdot\phi,\phi\right>\\ 
&=&  
\int_M| \nabla^{\theta}\phi|^2+ \dfrac{1}{4}\int_M \left({\rm R}|\phi|^2+
2\left<{\Theta},{\eta}^\phi\right>_0\right)\\
&\geq& 0, 
\end{eqnarray*}
so that $\phi$ is parallel, has non-zero constant length and no zeros.
Since
\[{\rm R}|\phi|^2+ 2\left<{\Theta},{\eta}^\phi\right>_0> 0\] 
at some point, 
\[0\geq  \int_M \left({\rm R}|\phi|^2+
2\left<{\Theta},{\eta}^\phi\right>_0\right)>0\]
which is a contradiction.\qd

{\bf Remark}.
This corollary means that one can check that a spinor is not harmonic by using the scalar
curvature of the manifold, the curvature operator of the connection on $F$ and the 2-forms
associated to the spinor.

\subsubsection{Killing spinors}
\begin{corol} 
Suppose $\phi\not=0$ is a Killing spinor with Killing constant $\mu$.
Then $\mu$ is either real or imaginary, and
\[
\mu^2
\quad\geq\quad {1\over 4n^2}\min_M({\rm R} -
2|\tilde\Theta^m|)
\quad\geq\quad {1\over 4n^2}\min_M({\rm R} -
2m|\tilde\Theta^1|). 
\]
If either of the two inequalities is attained, then $\phi$ is parallel, i.e. $\mu=0$. 
\end{corol}
{\em Proof}. Recall that
\begin{eqnarray*}
 \dirac^\theta (\phi) 
 &=& \sum_{i=1}^n e_i\cdot\nabla_{e_i}^\theta \phi \\
&=& -n\mu\, \phi .
\end{eqnarray*}
Then, by the twisted Schr\"odinger-Lichnerowicz formula \rf{SL-formula}
\[n^2\mu^2 \phi =   \Delta^{\theta}(\phi)+ \dfrac{\rm R}{4}\phi +
\dfrac{1}{2}\tilde\Theta^m\cdot\phi.\]
By taking  hermitian product with $\phi$ and integrating over $M$ we get
\begin{eqnarray*}
n^2\mu^2\int_M|\phi|^2
&=&  
\int_M| \nabla^\theta\phi|^2+ \int_M\dfrac{\rm  R}{4}|\phi|^2 +
\int_M\dfrac{1}{2}\left<\tilde\Theta^m\cdot\phi,\phi\right>\\
&\geq&\dfrac{1}{4}\int_M\left({\rm R} -
2|\tilde\Theta^m|\right)|\phi|^2\\
&\geq& {1\over4}\min_M({\rm R} -
2|\tilde\Theta^m|)\int_M|\phi|^2\\
&\geq& {1\over4}\min_M({\rm R} -
2m|\tilde\Theta^1|)\int_M|\phi|^2,
\end{eqnarray*}
and the inequalities follow. 
Since the right hand side of the equality above is a real number, $\mu$ must be either real or
imaginary.

Now, if either of the inequalities is attained,
\[\int_M|\nabla^\theta\phi|^2=0 \quad\quad\mbox{and}\quad\quad \nabla^\theta\phi=0.\]
\qd

\begin{corol}\label{Killing-lower-bound}
Suppose $\phi\not=0$ is a real Killing spinor with Killing constant $\mu$.
Then, 
\begin{eqnarray*}
\mu^2
&\geq& {1\over 4n^2{\rm vol}(M)}\int_M \left[{\rm R}+
{2\over |\phi|^2}\left<{\Theta},{\eta}^\phi\right>_0\right] ,
\end{eqnarray*}
and 
\begin{eqnarray*}
{\mu^2}
&\geq& {1\over 4n|\phi|^2{\rm vol}(M)}\int_M 2\left<{\Theta},{\eta}^\phi\right>_0 
-\left<\hat\Theta,\hat\eta^\phi\right>_1.
\end{eqnarray*}
\end{corol}
{\em Proof}. 
By the twisted Schr\"odinger-Lichnerowicz formula \rf{SL-formula}
\[n^2\mu^2 \phi =   \Delta^{\theta}(\phi)+ \dfrac{\rm R}{4}\phi +
\dfrac{1}{2}\tilde\Theta^m\cdot\phi.\]
By taking  hermitian product with $\phi$ and integrating we get
\begin{eqnarray*}
n^2\mu^2\int_M|\phi|^2
&=&  
\int_M\left< \Delta^{\theta}(\phi),\phi\right>+ \int_M\dfrac{\rm R}{4}\left<\phi,\phi\right> +
\dfrac{1}{2}\int_M\left<\tilde\Theta^m\cdot\phi,\phi\right>,\\ 
&\geq&  
 \dfrac{1}{4}\int_M {\rm R}|\phi|^2+
2\left<{\Theta},{\eta}^\phi\right>_0.
\end{eqnarray*}
Since $|\phi|$ is a non-zero constant
\begin{eqnarray*}
\mu^2
&\geq& {1\over 4n^2|\phi|^2{\rm vol}(M)}\int_M \left[{\rm  R}|\phi|^2+
2\left<{\Theta},{\eta}^\phi\right>_0\right] \\
&=& {1\over 4n^2|\phi|^2{\rm vol}(M)}\int_M \left[4n(n-1)\mu^2|\phi|^2
-\left<\hat\Theta,\hat\eta^\phi\right>_1
+2\left<{\Theta},{\eta}^\phi\right>_0\right] \\
&=& {(n-1)\mu^2\over n}
+{1\over 4n^2|\phi|^2{\rm vol}(M)}\int_M 2\left<{\Theta},{\eta}^\phi\right>_0 
-\left<\hat\Theta,\hat\eta^\phi\right>_1,
\end{eqnarray*}
where we have used Theorem \ref{Ricci-Killing-spinor}.
\qd

\subsubsection{Dirac eigen-spinors}

\begin{corol}
Suppose $\phi$ is
a Dirac eigenspinor
\[\dirac^\theta \phi = \lambda\phi.\]
for some $\lambda\in \mathbb{R}$.Then
\begin{eqnarray*}
\lambda^2
&\geq& {n\over4(n-1)}\left(\min_M({\rm R} -
2|\tilde\Theta^m|)\right)
\quad\geq\quad {n\over4(n-1)}\left(\min_M({\rm R} -
2m|\tilde\Theta^1|)\right).
\end{eqnarray*}
If either of the lower bounds is non-negative and is attained,
the spinor $\phi$ is a real Killing spinor with Killing constant
\[
\mu=\pm {1\over 2}\sqrt{{1\over n(n-1)}\min_M({\rm R} -
2|\tilde\Theta^m|)}\quad\quad\mbox{or}\quad\quad
\mu=\pm {1\over 2}\sqrt{{1\over n(n-1)}\min_M({\rm R} -
2m|\tilde\Theta^1|)},\]
respectively.
\end{corol}
{\em Proof}. 
Following \cite{friedrich}, let $h:M\lra\mathbb{R}$ be a fixed smooth function.
Consider the following metric connection on the twisted spin bundle
\[\nabla_X^h\phi = \nabla_X^\theta\phi+hX\cdot\phi.\]
Let
\[
\Delta^{h}(\phi) = - \sum_{i = 1}^n
\nabla^{h}_{e_i}\nabla^{h}_{e_i}\phi - \sum_{i = 1}
\diver(e_i)\nabla^{h}_{e_i}\phi, 
\]
be the Laplacian for this connection and recall that
\[|\nabla^h\phi|^2= \sum_{i=1}^n |\nabla_{e_i}^\theta\phi + he_i\cdot \phi|^2.\]
Then, by \rf{SL-formula}
\begin{eqnarray*}
 (\dirac^\theta - h)\circ (\dirac^\theta - h) (\phi)
&=&\dirac^\theta(\dirac^\theta\phi)-2h\dirac^\theta\phi-\grad(h)\cdot\phi+h^2\phi\\
&=&\Delta^{\theta}(\phi)+ \dfrac{\rm R}{4}\phi +
\dfrac{1}{2}\tilde\Theta^m\cdot\phi-2h\dirac^\theta\phi-\grad(h)\cdot\phi+h^2\phi.
\end{eqnarray*}
On the other hand,
\[\Delta^h\phi =  \Delta^\theta\phi -2h\dirac^\theta\phi-\grad(h)\cdot\phi+h^2\phi.\]
Thus
\[ 
(\dirac^\theta - h)\circ (\dirac^\theta - h) (\phi)
=\Delta^{h}(\phi)+ \dfrac{\rm R}{4}\phi +
\dfrac{1}{2}\tilde\Theta^m\cdot\phi+(1-n)h^2\phi 
\]
By using $\dirac^\theta\phi =\lambda\phi$, setting $h={\lambda\over n}$, taking hermitian product
with $\phi$ and integrating over $M$ we get
\[\lambda^2\left({n-1\over n}\right)^2\int_M|\phi|^2=\int_M|\nabla^{\lambda/n}\phi|^2+\lambda^2
{1-n\over n^2}\int_M|\phi|^2+\int_M\dfrac{\rm R}{4}|\phi|^2 +
\int_M\dfrac{1}{2}\left<\tilde\Theta^m\cdot\phi,\phi\right>\]
so that
\begin{eqnarray*}
\lambda^2\left({n-1\over n}\right)\int_M|\phi|^2
&=&\int_M|\nabla^{\lambda/n}\phi|^2+\int_M\dfrac{R}{4}|\phi|^2 +
\int_M\dfrac{1}{2}\left<\tilde\Theta^m\cdot\phi,\phi\right>\\ 
&\geq& {1\over4}\min_M({\rm R} -
2|\tilde\Theta^m|)\int_M|\phi|^2\\
&\geq& {1\over4}\min_M({\rm R} -
2m|\tilde\Theta^1|)\int_M|\phi|^2,
\end{eqnarray*}
and
\begin{eqnarray*}
\lambda^2
&\geq& {n\over4(n-1)}\min_M({\rm  R} -
2|\tilde\Theta^m|)
\quad\geq\quad {n\over4(n-1)}\min_M({\rm  R} -
2m|\tilde\Theta^1|).
\end{eqnarray*}

If either of the lower bounds is attained,
\[\int_M|\nabla^{\lambda/n}\phi|^2= 0,
\]
i.e.
\[\nabla^{\lambda/n}\phi= 0.
\]
\qd

Now, let $E\in\End(TM)$ be a fixed symmetric endomorphism and consider the following metric
connection on the twisted spin bundle
\[\nabla_X^E\phi = \nabla_X\phi+E(X)\cdot\phi.\]
Let
\[
\Delta^{E}(\phi) = - \sum_{i = 1}^n
\nabla^{E}_{e_i}\nabla^{E}_{e_i}(\phi) - \sum_{i = 1}
\diver(e_i)\nabla^{E}_{e_i}(\phi), 
\]
be this connection's Laplacian and 
\begin{eqnarray}
|\nabla^E\phi|^2
&=& \sum_{i=1}^n |\nabla_{e_i}\phi + E(e_i)\cdot \phi|^2\nonumber\\ 
&=& \sum_{i=1}^n |\nabla_{e_i}\phi |^2
-2{\rm Re}\left<E(e_i)\cdot\nabla_{e_i}\phi , \phi\right>
+|E(e_i)|^2 | \phi|^2\nonumber\\ 
&=&  |\nabla \phi |^2
+|E|^2 | \phi|^2 
-2{\rm Re}\sum_{i=1}^n\left<E(e_i)\cdot\nabla_{e_i}\phi , \phi\right>.\label{E-intermedia}
\end{eqnarray}
On the complement of the zero set of a spinor $\phi$, we can define a symmetric bilinear form 
$Q_\phi$ by
\[Q_\phi(X,Y) = {1\over 2}{\rm Re} \left<X\cdot \nabla_Y^\theta \phi + Y\cdot \nabla_X^\theta
\phi,{\phi\over |\phi|^2}\right>\] 
The associated field of
quadratic forms gives 
\[Q_\phi(e_i)= {\rm Re}\left<e_i\cdot\nabla_{e_i}^\theta \phi,{\phi\over |\phi|^2}\right>\]
and 
\[\tr(Q_\phi) = {\rm Re} \left<\dirac^\theta\phi,{\phi\over |\phi|^2}\right>\]
If $\phi\not=0$ is a Dirac spinor $\dirac^\theta(\phi) = \lambda\phi$, then
\[\tr(Q_\phi)=\lambda.\]
So, let us take
\begin{eqnarray*}
 E(X)
&=& (X\lrcorner Q_\phi)^\sharp\\
&=& \sum_{i=1}^nQ_\phi(X,e_i)e_i\\
&=& \ell^\phi(X),
\end{eqnarray*}
the so-called {\em energy-momentum tensor} of $\phi$.
Then, we can examine further the second and third summands of the identity \rf{E-intermedia} which
now looks as follows
\begin{eqnarray*}
|\nabla^{\ell^\phi}\phi|^2
&=&  |\nabla \phi |^2
+|\ell^\phi|^2 | \phi|^2 
-2{\rm Re}\sum_{i=1}^n\left<\ell^\phi(e_i)\cdot\nabla_{e_i}\phi , \phi\right>.
\end{eqnarray*}
On the one hand, 
\begin{eqnarray*}
|\ell^\phi|^2 
&=&  \sum_{i=1}^n |\ell^\phi(e_i)|^2 
\quad=\quad  \sum_{i,j=1}^n
Q_\phi(e_i,e_j)^2 ,
\end{eqnarray*}
and on the other,
\begin{eqnarray*}
-2{\rm Re}\sum_{i=1}^n\left<\ell^\phi(e_i)\cdot\nabla_{e_i}\phi , \phi\right> 
&=&  -2{\rm Re}\sum_{i=1}^n\left<\left(\sum_{j=1}^nQ_\phi(e_i,e_j)e_j\right)\cdot\nabla_{e_i}\phi ,
\phi\right> \\
&=&  -2|\phi|^2\sum_{i,j=1}^nQ_\phi(e_i,e_j)^2.
\end{eqnarray*}
Thus,
\begin{eqnarray*}
|\nabla\phi|^2
&=&  |\nabla^{\ell^\phi} \phi |^2
+|\ell^\phi|^2 | \phi|^2,
\end{eqnarray*}
so that
\begin{eqnarray*}
\int_M |\dirac^\theta \phi|^2
&=&  \int_M|\nabla^{\ell^\phi} \phi |^2
+|\ell^\phi|^2 | \phi|^2 
+ {{\rm R}\over 4}|\phi|^2 +{1\over 2}\left<\tilde\Theta^m\cdot\phi,\phi\right>.
\end{eqnarray*}
Since $\phi$ is a Dirac eigenspinor with eigenvalue $\lambda$,
\begin{eqnarray*}
\lambda^2\int_M|\phi|^2
&\geq&  \int_M
|\ell^\phi|^2 | \phi|^2 
+ {{\rm R}\over 4}|\phi|^2 +{1\over 2}\left<\tilde\Theta^m\cdot\phi,\phi\right>\\
&\geq& \min_{M}\left(|\ell^\phi|^2 + {{\rm R}\over 4} -{1\over 2}|\tilde\Theta^m|\right)
\int_M|\phi|^2\\
&\geq& \min_{M}\left(|\ell^\phi|^2 + {{\rm R}\over 4} -{m\over 2}|\tilde\Theta^1|\right)
\int_M|\phi|^2.
\end{eqnarray*}
i.e.
\begin{eqnarray*}
\lambda^2
&\geq& \min_{M}\left(|\ell^\phi|^2 + {R\over 4} -{1\over 2}|\tilde\Theta^m|\right) 
\quad\geq\quad \min_{M}\left(|\ell^\phi|^2 + {R\over 4} -{m\over 2}|\tilde\Theta^1|\right) .
\end{eqnarray*}

If the lower bound is attained, 
\[\int_M|\nabla^{\ell^\phi} \phi |^2=0\quad\quad\mbox{and}\quad\quad
\nabla^{\ell^\phi} \phi=0,\]
i.e.
\[\nabla_X\phi=-\ell^\phi(X)\cdot\phi\]
for all $X\in\Gamma(TM)$.
Furthermore, since $\nabla^{\ell^\phi}$ is compatible with the metric, $|\phi|$ is constant.

Thus, we have proved the following.

\begin{corol}
Suppose $\phi\not=0$ is
a Dirac eigenspinor
\[\dirac^\theta \phi = \lambda\phi.\]
Then
\begin{eqnarray*}
\lambda^2
&\geq& \min_{M}\left(|\ell^\phi|^2 + {R\over 4} -{1\over 2}|\tilde\Theta^m|\right) 
\quad\geq\quad \min_{M}\left(|\ell^\phi|^2 + {R\over 4} -{m\over 2}|\tilde\Theta^1|\right),
\end{eqnarray*}
where $\ell^\phi$ is the energy-momentum tensor of $\phi$.
If either of the lower bounds is non-negative and is attained,
$\phi$ has constant length and no zeros, and is a generalized Killing spinor with
symmetric endomorphism $\ell^\phi$, i.e.
\[\nabla_X\phi=-\ell^\phi(X)\cdot\phi\]
for all $X\in\Gamma(TM)$.
\end{corol}
\qd

{\small
\renewcommand{\baselinestretch}{0.5}
\newcommand{\bi}{\vspace{-.05in}\bibitem} }

\end{document}